\documentclass{article}
\date{}

\usepackage{amssymb}
\usepackage{amsfonts}
\newtheorem{theorem}{{\bf Theorem}}[section]
\newtheorem{lemma}[theorem]{{\bf Lemma}}
\newtheorem{corollary}[theorem]{{\bf Corollary}}
\newtheorem{remark}[theorem]{{\bf Remark}}

\newtheorem{question}[theorem]{{\bf Question}}

\newtheorem{definition}[theorem]{{\bf Definition}}
\newtheorem{proposition}[theorem]{{\bf Proposition}}
\newtheorem{example}[theorem]{{\bf Example}}

\newcommand{\proof}{{\bf Proof.  }}
\newcommand{\eproof}{{\square}}
\title{\Large \bf A few remarks on bounded operators \\on topological vector
spaces}

\author{{\bf Ljubi\v sa D.R. Ko\v cinac, Omid Zabeti}\\
University of Ni\v s, Faculty of Sciences and Mathematics,\\
18000 Ni\v s, Serbia\\
{\sf lkocinac@gmail.com}\\
Faculty of Mathematics, University of Sistan and Baluchestan,\\
P.O. Box 98135-674, Zahedan, Iran\\
{\sf o.zabeti@gmail.com}}

\begin{document}

\maketitle

\begin{abstract}
We give a few observations on different types of bounded operators
on a topological vector space $X$ and their relations with compact
operators on $X$. In particular, we investigate when these bounded
operators coincide with compact operators.  We also consider similar
types of bounded bilinear mappings between topological vector
spaces. Some properties of tensor product operators between locally
convex spaces are established. In the last part of the paper we deal
with operators on topological Riesz spaces.
\end{abstract}

\begin{flushleft}
\vspace{.3cm} {\sf Keywords}: Bounded operator; compact operator;
central operator; order bounded below operator; bounded bilinear
mapping; tensor product of locally convex spaces; Heine-Borel
property; topological Riesz space; locally solid Riesz space.
\end{flushleft}

%%%%%%%%%%%%%%%%%%%%%% 1111 %%%%%%%%%%%%
\section{Introduction}
%%%%%%%%%%%%%%%%% 1111 %%%%%%%%%%%%%%%%%

Throughout the paper, all topological vector spaces are over the
scalar field $\mathbb K$ which is either the field $\mathbb R$ of
real numbers or the field $\mathbb C$ of complex numbers. A
neighborhood $U$ of $0$ in a topological vector space will be simply
called a zero neighborhood.

Let $X$ be a topological vector space. A subset $A$ of $X$ is
\emph{bounded} if for each zero neighborhood $U$ in $X$ there is a
scalar $\lambda$ such that $A \subseteq \lambda U$. $X$ is said to
be \emph{locally bounded} if there is a bounded neighborhood of
$0\in X$. In the literature one can find two different notions of
bounded operators. In \cite{KN, Ru, Sch} there is the definition
which is a definition of {\sf nb}-bounded operators below, while in
\cite{Edw} one can find the definition which is a definition of {\sf
bb}-bounded operators; the following terminology is from \cite{Tr}.
A linear operator $T$ on $X$ is said to be:

\begin{itemize}
\item \emph{${\sf nb}$-bounded} if there exists a zero neighborhood
$U\subseteq X$ such that $T(U)$ is a bounded subset of $X$;

\item \emph{${\sf bb}$-bounded} if $T$ mappings bounded sets into
bounded sets.
\end{itemize}

The class of all ${\sf nb}$-bounded operators on $X$ is denoted by
${\sf B_n}(X)$ and is equipped with the topology of uniform
convergence on some zero neighborhood. A net $(S_{\alpha})$ of ${\sf
nb}$-bounded operators is said to converge to zero uniformly on a
zero neighborhood $U \subseteq X$ if for each zero neighborhood $V$,
there is an $\alpha_0$ with $S_{\alpha}(U) \subseteq V$ for each
$\alpha\ge\alpha_0$.

The class of all ${\sf bb}$-bounded operators on $X$ is denoted by
${\sf B_b}(X)$ and is endowed with the topology of uniform
convergence on bounded sets. Note that a net $(S_{\alpha})$ of ${\sf
bb}$-bounded operators converges to zero uniformly on a bounded set
$B \subseteq X$ if for each zero neighborhood $V$, there is an index
$\alpha_0$ with $S_{\alpha}(B) \subseteq V$ for all
$\alpha\ge\alpha_0$.

Also, the class of all continuous operators on $X$ is denoted by
${\sf B_c}(X)$ and is assigned to the topology of equicontinuous
convergence. A net$(S_{\alpha})$ of continuous operators is
convergent to the zero operator equicontinuously if for each zero
neighborhood $V\subseteq X$ there exists a zero neighborhood $U$
such that for every $\varepsilon>0$ there is an $\alpha_0$ with the
property $S_{\alpha}(U)\subseteq \varepsilon V$ for any
$\alpha\ge\alpha_0$.

It is easy to see that
\[
{\sf B_n}(X)\subseteq {\sf B_c}(X)\subseteq {\sf B_b}(X).
\]
Also note that the above inclusions become equalities when $X$ is
locally bounded; see \cite{Tr, Omid} for more details concerning
these operators. Recall that a topological vector space $X$ has the
\emph{Heine-Borel property} if every closed and bounded subset of
$X$ is compact.

By $X\otimes Y$, we mean the algebraic tensor product space. If $X$
and $Y$ are locally convex spaces, then the symbol
$X{\otimes}_{\pi}Y$ will be used for the algebraic tensor product
space endowed with the projective tensor product topology. For a
review about projective tensor product of locally convex spaces and
the related notions, we refer the reader to \cite[Chapter III; 5,
6]{Sch}; for a topological flavour of topological vector spaces and
the relevant aspects, one may consult \cite[Chapters I, II,
III]{Ru}. Note that the symbol ${\rm co}(A\otimes B)$ denotes the
\emph{convex hull} of $A\otimes B$.

Some notions will be given in the course of exposition in the
beginning of a section.

The paper is organized in the following way. In Section 2 we further
investigate bounded linear operators on topological vector spaces.
In Section 3 we introduce bounded bilinear mappings between
topological vector spaces and consider some their properties.
Section 4 is devoted to a topological approach to the notions of
central and order bounded below operators defined on a topological
Riesz space. With an appropriate topology, we extend some known
results for central and order bounded below operators on a Banach
lattice, to the topological Riesz space setting.

%%%%%%%%%%%%%%%% 2222 %%%%%%%%%%%%

\section{Linear mappings}

%%%%%%%%%%%%%%%% 2222 %%%%%%%%%%%%

Let $X$ be a normed space, ${\sf K}(X)$ be the space of all compact
operators on $X$, and ${\sf B}(X)$ be the collection of all bounded
linear operators. From the equality $\|T(x)\|\leq \|T\|\|x\|$, it
follows that $X$ is a topological ${\sf B}(X)$-module, where the
module multiplication is given via the formula $(T,x)\mapsto T(x)$,
for each linear operator $T$ and each $x\in X$. On the other hand,
it is known that ${\sf K}(X)$ is a closed subspace of ${\sf B}(X)$.
It is also of interest to investigate situations in which ${\sf
K}(X)$ and ${\sf B}(X)$ are the same. So, it is natural to see if
these results can be generalized to ordinary topological vector
spaces and to different classes of bounded and compact operators on
them.

In \cite{Tr}, two different notions for compact operators on a
topological vector space have been introduced. A linear operator $T$
on a topological vector space $X$ is said to be:
\begin{itemize}
\item \emph{${\sf n}$-compact} if there is a zero neighborhood
$U\subseteq X$ for which $T(U)$ is relatively compact (which means
that its closure is compact);

\item \emph{${\sf b}$-compact} if for
each bounded set $B\subseteq X$, $T(B)$ is relatively compact.
\end{itemize}
We use the notations ${\sf K_n}(X)$ and ${\sf K_b}(X)$ for the set
of all ${\sf n}$-compact linear operators and the set of all ${\sf
b}$-compact linear operators on $X$, respectively. It is easy to see
that ${\sf K_n}(X)$ is a two-sided ideal of ${\sf B_n}(X)$ and ${\sf
K_b}(X)$ is a right ideal of ${\sf B_b}(X)$.

In this section we investigate some relations between bounded linear
operators and compact ones. For more details about bounded and
compact operators on topological vector spaces and the related
notions see \cite{Edw, Sch, Tr, Omid}.

In \cite{Omid}, it has been proved that each class of bounded linear
operators on a topological vector space $X$, with respect to the
appropriate topology, forms a topological algebra. In this section
we show that $X$ is a topological $A$-module, where $A$ is one of
the topological algebras ${\sf B_n}(X)$, ${\sf B_c}(X)$, and ${\sf
B_b}(X)$, respectively, and the module multiplication is given via
$(T,x)\to T(x)$, for every linear operator $T$ and every $x\in X$.
Let us point out that a topological vector space $X$ is a
\emph{topological $A$-module}, where $A$ is a topological algebra
over the same field, provided that the module multiplication is
continuous as a mapping from $A\times X$, equipped with product
topology, into $X$.

\begin{proposition}
The module multiplication in ${\sf B_n}(X)$ is continuous with
respect to the topology of uniform convergence on some zero
neighborhood.
\end{proposition}
$\proof$ Let $(x_{\alpha})$ be a net in $X$ which is convergent to
zero and $(T_{\alpha})$ be a net in ${\sf B_n}(X)$ converging to
zero uniformly on some zero neighborhood $U\subseteq X$. There is an
$\alpha_0$ such that $x_\alpha \in U$ for each $\alpha\geq\alpha_0$.
If $V$ is an arbitrary zero neighborhood in $X$, Then, there exists
an index $\alpha_1$ with $T_{\alpha}(U)\subseteq V$ for each
$\alpha\geq\alpha_1$, so that for sufficiently large $\alpha$, we
have
\[
T_{\alpha}(x_{\alpha})\subseteq T_{\alpha}(U)\subseteq V.
\]
$\eproof$

\begin{proposition}
The module multiplication in ${\sf B_c}(X)$ is continuous with
respect to the equicontinuous convergence topology.
\end{proposition}
$\proof$ Let $(x_{\alpha})$ be a net in $X$ which is convergent to
zero and $(T_{\alpha})$ be a net of continuous operators converging
to zero equicontinuously. Suppose $V$ is an arbitrary zero
neighborhood in $X$. There exist a zero neighborhood $U\subseteq X$
and an index $\alpha_0$ with $T_{\alpha}(U)\subseteq V$ for each
$\alpha\geq\alpha_0$. Choose an $\alpha_1$ such that $x_\alpha \in
U$ for every $\alpha\geq\alpha_1$. Thus, for sufficiently large
$\alpha$, we conclude
\[
T_{\alpha}(x_{\alpha})\subseteq T_{\alpha}(U)\subseteq V.
\]
$\eproof$

\begin{proposition}
The module multiplication in ${\sf B_b}(X)$ is sequentially
continuous with respect to the topology of uniform convergence on
bounded sets.
\end{proposition}
$\proof$ Let $(x_n)$ be a sequence in $X$ which is convergent to
zero and $(T_n)$ be a sequence of ${\sf bb}$-bounded operators
converging to zero uniformly on bounded sets. Note that
$E=\{x_n:n\in \mathbb N\}$ is bounded in $X$. Suppose $V$ is an
arbitrary zero neighborhood in $X$. There exists an $n_0$ such that
$T_n(E)\subseteq V$ for any $n>n_0$, so that we have
\[
T_n(x_n)\subseteq T_n(E)\subseteq V.
\]
$\eproof$

\begin{question} Is the module multiplication in ${\sf B_b}(X)$
continuous, in general?
\end{question}

Note that for a normed space $X$, ${\sf K}(X)={\sf B}(X)$ if and
only if $X$ is finite dimensional. In this step, we consider some
situations where a class of compact linear operators coincides with
the corresponding class of bounded operators.

\begin{proposition}\label{3}
${\sf K_{b}}(X)={\sf B_{b}}(X)$ if and only if $X$ has the
Heine-Borel property.
\end{proposition}
$\proof$ ${\sf K_{b}}(X)={\sf B_{b}}(X)$ $\Leftrightarrow$ $I \in
{\sf K_{b}}(X)$ $\Leftrightarrow$ $X$ has the Heine-Borel property,
where $I$ denotes the identity operator on $X$. $\eproof$

\begin{remark}\label{4}\rm
Note that when $X$ has the Heine-Borel property, then ${\sf
K_{n}}(X)={\sf B_{n}}(X)$. Suppose, for a topological vector space
$X$, ${\sf K_{n}}(X)={\sf B_{n}}(X)$. We consider two cases. First,
assume $X$ is locally bounded. Then,

$I\in {\sf B_n}(X)$ $\Rightarrow$ $I\in {\sf K_n}(X)$ $\Rightarrow$
$X$ is locally compact $\Rightarrow$ $X$ is finite dimensional.

The second case, when $X$ is not locally bounded. Then,

$I\not\in {\sf B_n}(X)$ $\Rightarrow$ $I\not\in {\sf K_n}(X)$
$\Rightarrow$ $X$ is not locally compact $\Rightarrow$ $X$ is
infinite dimensional.
\end{remark}

Let $(X_n)$ be a sequence of topological vector spaces in which,
every $X_n$ has the Heine-Borel property. Put
$X=\prod_{n=1}^{\infty}X_n$, with the product topology. It is known
that $X$ is a topological vector space. In the following, we
establish that each $X_n$ has the Heine-Borel property if and only
if so has $X$.

\begin{theorem}\label{1}
Let $X=\prod_{n=1}^{\infty}X_n$, with the product topology; then $X$
has the Heine-Borel property if and only if each $X_n$ has this
property, as well.
\end{theorem}
$\proof$ First, assume that each $X_n$ has the Heine-Borel property.
We claim that if $B\subseteq X$ is bounded, then there exist bounded
subsets $B_i\subseteq X_i$ such that $B\subseteq
\prod_{i=1}^{\infty} B_i$. Put
 \[
 B_i=\{x\in X_i: (y_1,\ldots,y_{i-1},x,y_{i+1},\ldots)\in B, y_j\in X_j\}.
 \]
Each $B_i$ is bounded in $X_i$. Let $W_i$ be a zero neighborhood in
$X_i$ and put
 \[
 W=X_1\times\ldots X_{i-1}\times W_i\times X_{i+1}\times\ldots.
 \]
Since $W$ is a zero neighborhood in $X$, there exists a positive
scalar $\alpha$ such that $B\subseteq \alpha W$, so that
$B_i\subseteq \alpha W_i$. Also, it is easy to see that $B\subseteq
\prod_{i=1}^{\infty} B_i$. Therefore,
 \[
 \overline{B}\subseteq \overline{\prod_{i=1}^{\infty} B_i}= \prod_{i=1}^{\infty}\overline{B_i},
 \]
so that we conclude $\overline{B}$ is also compact, i.e. that $X$
has the Heine-Borel property.

For the converse, suppose $X$ has the Heine-Borel property. Choose a
bounded set $B_n\subseteq X_n$. Put
 \[
 B=\{0\}\times\ldots\times\{0\}\times B_n\times\{0\}\times\ldots.
 \]
It is an easy job to see that $B$ is bounded in $X$, so that
$\overline{B}$ is compact. By using Tychonoff's theorem, we conclude
that $\overline{B_n}$ is compact and this would complete our claim.
$\eproof$

\medskip
Collecting results of Theorem \ref{1}, Proposition \ref{3}, and
Remark \ref{4}, we have the following.

\begin{corollary}
Let $(X_n)$ be a sequence of topological vector spaces, in which,
each $X_n$ has the Heine-Borel property. Put
$X=\prod_{n=1}^{\infty}X_n$, with the product topology. Then, ${\sf
K_b}(X)={\sf B_b}(X)$ and ${\sf K_n}(X)={\sf B_n}(X)$.
\end{corollary}

\begin{remark}\rm
Note that when $X$ has the Heine-Borel property, ${\sf K_n}(X)$ need
not be equal to ${\sf K_b}(X)$. For example, consider the identity
operator on ${\mathbb R}^{\mathbb N}$. Indeed, it is ${\sf
b}$-compact but it fails to be $n$-compact; nevertheless, $X$ has
the Heine-Borel property.
\end{remark}

Compact operators are not closed in the  topologies induced by the
corresponding class of bounded linear operators. To see this,
consider the following examples.

\begin{example}\rm
${\sf K_{n}}(X)$ is not a closed subspace of ${\sf B_{n}}(X)$, in
general. Let $X$ be $c_{00}$, the space of all real null sequences,
with the uniform norm topology. Suppose that $T_n$ is the linear
operator defined by
\[
T_n(x_1,x_2,\ldots,x_n,\ldots)=(x_1,\frac{1}{2}x_2,\ldots,
\frac{1}{n}x_n,0,\ldots).
\]
It is easy to see that each $T_n$ is ${\sf n}$-compact. Also,
$(T_n)$ converges uniformly on $N_{1}^{(0)}$, the open unit ball of
$X$ with center zero, to the linear operator $T$ defined by
$T(x_1,x_2,\ldots)=(x_1,\frac{1}{2}x_2,\ldots)$. For, if
$\varepsilon>0$ is arbitrary, there is an $n_0 \in \mathbb N$ such
that $\frac{1}{n_0}<\varepsilon $. So, for each $n>n_0$,
$(T_n-T)(N_{1}^{(0)})\subseteq N_{\varepsilon}^{(0)}$. Now, it is
not difficult to see that $T$ is ${\sf nb}$-bounded but it is not an
${\sf n}$-compact linear operator.
\end{example}

\begin{example}\rm
${\sf K_{b}}(X)$ fails to be closed in ${\sf B_{b}}(X)$, in general.
Let $X$ be $c_0$, the space of all vanishing sequences, with the
coordinate-wise topology. Let $P_n$ be the projection on the first
$n$-components. By using the Tychonoff's theorem, we can conclude
that for each $n \in \mathbb N$, $P_n$ is ${\sf b}$-compact. Also,
$(P_n)$ converges uniformly on bounded sets to the identity operator
$I$. We show that $I$ is not ${\sf b}$-compact. Suppose that $B$ is
the sequence $(a_n)$
 defined by $a_n=(1,1,\ldots,1,0,0,\ldots)$ in which 1 is appeared $n$ times.
 $B$ is a Cauchy sequence in $c_0$, so that it is bounded. Also note that
 $\overline{B}=B$. Now, if $I \in {\sf K_{b}}(c_0)$, then $B$
should be compact. Since $B$ is not complete, this is impossible.
\end{example}

%%%%%%%%%%%%%%%%% 3333 %%%%%%%
\section{Bilinear mappings}

%%%%%%%%%%%%%%%%% 3333 %%%%%%%
The notion of a jointly continuous bilinear mapping between
topological vector spaces has been studied widely, for example, see
\cite{Ru, Ryan, Sch} for more information. In particular, when we
deal with the normed spaces framework, these mappings carry bounded
sets (with respect to the product topology) to bounded sets. On the
other hand, tensor products are a fruitful and handy tool in
converting a bilinear mapping to a linear operator in any setting;
for example, the projective tensor product for normed spaces and the
Fremlin projective tensor products for vector lattices and Banach
lattices (see \cite{Fremlin1, Fremlin2} for ample information). In a
topological vector space setting, we can consider two different
non-equivalent ways to define a bounded bilinear mapping. It turns
out that these aspects of boundedness are in a sense "intermediate"
notions of a jointly continuous one. On the other hand, different
types of bounded linear operators between topological vector spaces
and some of their properties have been investigated (see \cite {Tr,
Omid}). In this section, by using the concept of projective tensor
product between locally convex spaces, we show that, in a sense,
different notions of a bounded bilinear mapping coincide with
different aspects of a bounded operator. We prove  that for two
bounded linear operators, the tensor product operator also has the
same boundedness property, as well.

\begin{definition}\rm
Let $X$, $Y$, and $Z$ be topological vector spaces. A bilinear
mapping $\sigma:X \times Y \to Z$ is said to be:
\begin{itemize}
\item[$(i)$] \emph{${\sf n}$-bounded} if there exist some zero neighborhoods $U\subseteq X$
and $V\subseteq Y$ such that  $\sigma(U\times V)$ is bounded in $Z$;
\item[$(ii)$]\emph{${\sf b}$-bounded} if for any bounded sets $ B_1\subseteq X$ and $B_2\subseteq Y$,
$\sigma(B_1\times B_2)$ is bounded in $Z$.
\end{itemize}
\end{definition}

We first show that these concepts of bounded bilinear mappings are
not equivalent.

\begin{example}\rm
Let $X = {\mathbb R}^{\mathbb N}$ be the space of all real sequences
with the Tychonoff product topology. Consider the bilinear mapping
$\sigma: X\times X\to X$ defined by $\sigma(x,y)=xy$ where
$x=(x_i)$, $y=(y_i)$ and the product is pointwise. It is easily
verified that $\sigma$ is ${\sf b}$-bounded; but since $X$ is not
locally bounded, it can not be an ${\sf n}$-bounded bilinear
mapping.
\end{example}

It is not difficult to see that every ${\sf n}$-bounded bilinear
mapping is jointly continuous and every jointly continuous bilinear
mapping is ${\sf b}$-bounded, so that these concepts of bounded
bilinear mappings are related to jointly continuous bilinear
mappings. Note that a ${\sf b}$-bounded bilinear mapping need not be
jointly continuous, even separately continuous; by a separately
continuous bilinear mapping, we mean one which is continuous in each
of its components. Consider the following example (which is actually
originally an exercise from \cite[Chapter I, Exercise 13]{Ru}; we
will give a proof for it for the sake of completeness).

\begin{example} \rm
Let $X$ be the space $C[0,1]$, consisting of all real continuous
functions on $[0,1]$. Suppose $\tau_1$ is the topology generated by
the seminorms $p_{x}(f)=|f(x)|$, for each $x\in [0,1]$, and $\tau_2$
is the topology induced by the metric defined via the formulae
\[
d(f,g)=\int_{0}^{1}\frac{|f(x)-g(x)|}{1+|f(x)-g(x)|}dx.
\]
Consider the bilinear mapping $\sigma: (X,\tau_1)\times
(X,\tau_1)\to (X,\tau_2)$ defined by $\sigma(f,g)=fg$. It is easy to
show that $\sigma$ is a ${\sf b}$-bounded bilinear mapping. But is
it not even separately continuous; for example the mapping $g=1_X$,
the identity operator from $(X,\tau_1)$ into $(X,\tau_2)$, is not
continuous. To see this, suppose
\[
V=\{f\in X: d(f,0)<\frac{1}{2}\}.
\]
$V$ is a zero neighborhood in $(X,\tau_2)$. If the identity operator
is continuous, there should be a zero neighborhood $U\subseteq
(X,\tau_1)$ with $U\subseteq V$. Therefore, there are
$\{x_1,\ldots,x_n\}$ in $[0,1]$ and $\varepsilon>0$ such that
\[
U=\{f\in X, |f(x_i)|< \varepsilon,  i=1,\ldots,n\}.
\]
For each subinterval $[x_i,x_{i+1}]$, consider positive reals
$\alpha_i$ and $\alpha_{i+1}$ such that $x_i<\alpha_i<
\alpha_{i+1}< x_{i+1}$. For an $n\in \mathbb N$, Define,
\[
f_i(x) =    \left \{\begin{array}{rl}
\frac{n(x-x_i)}{\alpha_i-x_i}, & \mbox{ if } x_i \leq x \leq
\alpha_i,\\ n, & \mbox{ if } \alpha_i
\leq x \leq \alpha_{i+1},\\
\frac{n(x-x_{i+1})}{\alpha_{i+1}-x_{i+1}}, & \mbox{ if }
\alpha_{i+1} \leq x \leq x_{i+1}. \end{array}
        \right .
\]
Now consider the continuous function $f$ on $[0,1]$ defined by
${f_i}'s$. Obviously, $f\in U$. Put
$\beta=\sum_{i=1}^{n}(\alpha_{i+1}-\alpha_i)$. We can choose $n\in
\mathbb N$ and $\beta$ in such a way that $\frac{\beta
n}{n+1}>\frac{1}{2}$. Thus,
\[
\int_{0}^{1}\frac{|f(x)|}{1+|f(x)|}dx> \frac{\beta
n}{n+1}>\frac{1}{2}.
\]
This completes the claim.
\end{example}

In what follows, by using the concept of the projective tensor
product of locally convex spaces, we are going to show that these
concepts of bounded bilinear mappings are, in fact, the different
types of bounded operators defined on a locally convex topological
vector space. Recall that if $U$ and $V$ are zero neighborhoods for
locally convex spaces $X$ and $Y$, respectively, then ${\rm
co}(U\otimes V)$ is a typical zero neighborhood for the locally
convex space $X{\otimes}_{\pi}Y$.

\begin{proposition}
Let $X$, $Y$ and $Z$ be locally convex vector spaces and $\theta:
X\times Y\to X{\otimes}_{\pi} Y$ be the canonical bilinear mapping.
If $\varphi: X\times Y\to Z$ is an ${\sf n}$-bounded bilinear
mapping, there exists an ${\sf nb}$-bounded operator
$T:X{\otimes}_{\pi} Y\to Z$ such that $T\circ \theta=\varphi$.
\end{proposition}
$\proof$ By \cite[III.6.1]{Sch}, there is a linear mapping
$T:X{\otimes}_{\pi} Y\to Z$ such that $T\circ \theta=\varphi$.
Therefore, it is enough to show that $T$ is ${\sf nb}$-bounded.
Since $\varphi$ is ${\sf n}$-bounded, there are zero neighborhoods
$U\subseteq X$ and $V\subseteq Y$ such that $\varphi(U\times V)$ is
bounded in $Z$. Let $W$ be an arbitrary zero neighborhood in $Z$.
There is $r>0$ with $\varphi(U\times V)\subseteq rW$. It is not hard
to show that $T(U\otimes V)\subseteq rW$, so that $T({\rm
co}(U\otimes V))\subseteq rW$. But, by the fact mentioned before
this proposition, ${\rm co}(U\otimes V)$ is a zero neighborhood in
the space $X{\otimes}_{\pi}Y$. This completes the proof. $\eproof$

\begin{proposition}\label{1}
Let $X$, $Y$ and $Z$ be locally convex vector spaces and $\theta:
X\times Y\to X{\otimes}_{\pi} Y$ be the canonical bilinear mapping.
If $\varphi: X\times Y\to Z$ is a ${\sf b}$-bounded bilinear
mapping, there exists a ${\sf bb}$-bounded operator
$T:X{\otimes}_{\pi} Y\to Z$ such that $T \circ\theta=\varphi$.
\end{proposition}
$\proof$ As in the proof of the previous theorem, the existence of
the linear mapping $T:X{\otimes}_{\pi} Y\to Z$ such that $T
\circ\theta=\varphi$ follows by \cite[III.6.1]{Sch}. We prove that
the linear mapping $T$ is ${\sf bb}$-bounded. Consider a  bounded
set $B\subseteq X{\otimes}_{\pi} Y$. There exist bounded sets
$B_1\subseteq X$ and $B_2\subseteq Y$ such that $B\subseteq
B_1\otimes B_2$. To see this, put
\begin{center}
$B_1=\{x\in X, \exists\ y\in Y, $ {such that} $x\otimes y\in B\},$
\end{center}
\begin{center}
$B_2=\{y\in Y, \exists\ x\in X, $ {such that} $ x\otimes y\in B\}.$
\end{center}
It is not difficult to see that $B_1$ and $B_2$ are bounded in $X$
and $Y$, respectively, and $B\subseteq B_1\otimes B_2$. Also, since
$\theta$ is jointly continuous, $B_1\otimes B_2$ is also bounded in
$X{\otimes}_{\pi}Y$. Thus, from the inclusion
\[
T(B)\subseteq T(B_1\otimes B_2)=T \circ\theta(B_1\times
B_2)=\varphi(B_1\times B_2)
\]
and using the fact that $\varphi$ is a ${\sf b}$-bounded bilinear
mapping, it follows that $T$ is a ${\sf bb}$-bounded linear
operator. This concludes the claim and completes the proof of the
proposition. $\eproof$

\begin{remark} \rm
Note that the similar result for jointly continuous bilinear
mappings between locally convex spaces is known and commonly can be
found in the contexts concerning topological vector spaces (see for
example \cite[III.6]{Sch}).
\end{remark}

We are going now to investigate whether or not the tensor product of
two operators preserves different kinds of bounded operators between
topological vector spaces. The response is affirmative. Recall that
for vector spaces $X$, $Y$, $Z$, and $W$, and linear operators
$T:X\to Y$,  $S:Z\to W$, by the tensor product of $T$ and $S$, we
mean the unique linear operator $T\otimes S: X\otimes Z\to Y\otimes
W$ defined via the formulae
\[
(T\otimes S)(x\otimes z)=T(x)\otimes S(z);
\]
one may consult \cite{Ryan} for a comprehensive study regarding the
tensor product operators.

\begin{theorem}
Let $X$, $Y$, $Z$, and $W$ be locally convex spaces, and $T:X\to Y$
and $S:Z\to W$ be ${\sf nb}$-bounded linear operators. Then the
tensor product operator $T\otimes S: X{\otimes}_{\pi} Z\to
{Y\otimes}_{\pi} W$ is ${\sf nb}$-bounded.
\end{theorem}
$\proof$ Let $U\subseteq X$ and $V\subseteq Z$ be two zero
neighborhoods such that $T(U)$ and $S(V)$ are bounded subsets of $Y$
and $W$, respectively. Let $O_1\subseteq Y$ and $O_2\subseteq W$ be
two arbitrary zero neighborhoods. There exist positive reals
$\alpha$ and $\beta$ with $T(U)\subseteq \alpha O_1$ and
$S(V)\subseteq \beta O_2$. Then
\[
(T\otimes S)(U\otimes V)=T(U)\otimes S(V)\subseteq \alpha\beta
(O_1\otimes O_2)\subseteq \alpha\beta {\rm co}(O_1\otimes O_2),
\]
so that $(T\otimes S)({\rm co}(U\otimes V))\subseteq \alpha\beta
{\rm co}(O_1\otimes O_2)$. This is the desired result. $\eproof$

\begin{theorem}
Suppose $X$, $Y$, $Z$, and $W$ are locally convex spaces, and
$T:X\to Y$ and $S:Z\to W$ are ${\sf bb}$-bounded linear operators.
Then the tensor product operator $T\otimes S: X{\otimes}_{\pi} Z\to
Y{\otimes}_{\pi} W$ is also ${\sf bb}$-bounded.
\end{theorem}
$\proof$ Fix a bounded set $B\subseteq X{\otimes}_{\pi} Z$. By the
argument used in Proposition \ref{1}, there are bounded sets
$B_1\subseteq X$ and $B_2\subseteq Z$ with $B\subseteq B_1\otimes
B_2$. Let $O_1\subseteq Y$ and $O_2\subseteq W$ be two arbitrary
zero neighborhoods. There are positive reals $\alpha$ and $\beta$
such that $T(B_1)\subseteq \alpha O_1$ and $S(B_2)\subseteq \beta
O_2$. Therefore,
\[
(T\otimes S)(B)\subseteq (T\otimes S)(B_1\otimes B_2)=T(B_1)\otimes
S(B_2)\subseteq \alpha\beta(O_1\otimes O_2)\subseteq \alpha\beta
{\rm co}(O_1\otimes O_2),
\]
hence $(T\otimes S)(B)\subseteq \alpha\beta{\rm co}(O_1\otimes
O_2)$, as required. $\eproof$

\begin{theorem}
Suppose $X$, $Y$, $Z$, and $W$ are locally convex spaces, and
$T:X\to Y$ and $S:Z\to W$ are continuous linear operators. Then the
tensor product operator $T\otimes S: X{\otimes}_{\pi} Z\to
{Y\otimes}_{\pi} W$ is jointly continuous.
\end{theorem}
$\proof$ Let $O_1\subseteq Y$ and $O_2\subseteq Z$ be two arbitrary
zero neighborhoods. There exist zero neighborhoods $U\subseteq X$
and $V\subseteq Z$ such that $T(U)\subseteq  O_1$ and $S(V)\subseteq
O_2$. It follows
\[
(T\otimes S)(U\otimes V)=T(U)\otimes S(V)\subseteq (O_1\otimes
O_2)\subseteq {\rm co}(O_1\otimes O_2),
\]
so that $(T\otimes S)({\rm co}(U\otimes V))\subseteq {\rm
co}(O_1\otimes O_2)$, as claimed. $\eproof$

%%%%%%%%%%%%%%%%%% 44444 %%%%%%%%%%%%%%%%%%%%%%%%%%%%
\section{Operators in topological Riesz spaces}
%%%%%%%%%%%%%%%%%% 44444 %%%%%%%%%%%%%%%%%%%%%%%%%%%%

In this section we give a topological approach to the notions of
central operators and order bounded below operators defined on a
topological Riesz space. With an appropriate topology, we extend to
topological Riesz spaces some known results for these operators on
Banach lattices.

We recall some concepts and terminology. A \emph{Riesz space} (or
\emph{vector lattice}) is an ordered real vector space $X$ which is
also a lattice. For $x\in X$ one defines $x^+ = x\vee 0$ (the
\emph{positive part} of $x$), $x^- = (-x)\vee 0$ (the \emph{negative
part} of $x$), and $|x| = x \vee (-x) = x^+ + x^-$ (the
\emph{absolute value} or \emph{modulus} of $x$).

A subset $S$ of a Riesz space $X$ is said to be \textit{solid} if
$y\in S$, $x\in X$, and $|x|\leq |y|$ imply $x\in S$. A
\textit{topological Riesz space} is a Riesz space which is at the
same time a (ordered) topological vector space. By a \textit{locally
solid Riesz space} we mean a topological Riesz space with a locally
solid topology. A \textit{Banach lattice} is a Riesz space which is
also a Banach space, where the norm is a lattice norm.

For more information on topological Riesz spaces and related
notions, and also for terminology used in this section, we refer the
reader to \cite{AA, AAP, AB, AB1, fremlin}.

\medskip
\subsection{Central operators}

\medskip
We recall that a linear operator $T$ on a Riesz space $X$ is called
\emph{central} if there exists a positive real number $\gamma$ such
that for each $x\in X$, we have $|T(x)|\leq \gamma|x|$ (see
\cite{BK, CD, AW}, where there are also interesting results
concerning this class of operators). In \cite{AW}, Wickstead showed
that if $X$ is a Banach lattice, then ${\sf Z}(X)$, the space of all
central operators on $X$ (called the \emph{center} of $X$), is a
unital Banach algebra with respect to the norm given by
\[
\|T\|=\inf\{\lambda\ge 0: |T|\leq \lambda I\},
\]
for each $T\in {\sf Z}(X)$, where $I$ denotes the identity operator
on $X$.

We are going to generalize this result to central operators on a
locally solid Riesz space $X$ endowed with the
\emph{$\tau$-topology}: a net $(T_{\alpha})$ of central operators
converges to zero in the $\tau$-topology if for each $\varepsilon>0$
there is an $\alpha_0$ such that $|T_{\alpha}(x)|\leq \varepsilon
|x|$, for each $\alpha\ge\alpha_0$ and for each $x\in X$. It is easy
to see that ${\sf Z}(X)$ is a unital algebra. We show that ${\sf
Z}(X)$ is in fact a topological algebra.

\begin{proposition}\label{Z(X)1}
The operations of addition, scalar multiplication and product are
continuous in ${\sf Z}(X)$ with respect to the $\tau$-topology.
\end{proposition}
$\proof$ Suppose $(T_{\alpha})$ and $(S_{\alpha})$ are two nets of
central operators which are convergent to zero in the
$\tau$-topology. Let $\varepsilon>0$ be given. There are indices
$\alpha_0$ and $\alpha_1$ such that $|T_{\alpha}(x)|\leq
\frac{\varepsilon}{2}|x|$ for each $\alpha\geq\alpha_0$ and $x\in
X$, and $|S_{\alpha}(x)|\leq \frac{\varepsilon}{2}|x|$ for each
$\alpha\geq\alpha_1$ and $x\in X$. Choose $\alpha_2$ with
$\alpha_2\geq\alpha_0$ and $\alpha_2\geq\alpha_1$. Then for each
$\alpha\geq\alpha_2$ and for each $x\in X$, we have,
\[
|(T_{\alpha}+S_{\alpha})(x)|\leq |T_{\alpha}(x)|+|S_{\alpha}(x)|\leq
\varepsilon|x|.
\]
Now, we show the continuity of the scalar multiplication. Suppose
$(\gamma_{\alpha})$ is a net of reals which converges to zero.
Without loss of generality, we may assume that
$|\gamma_{\alpha}|\leq 1$ for each $\alpha$. Therefore, for each
$x\in X$ we have
\[
|\gamma_{\alpha}||T_{\alpha}(x)|\leq
|\gamma_{\alpha}|\varepsilon|x|\leq \varepsilon|x|,
\]
for all $\alpha\ge\alpha_0$.

\noindent For continuity of the product, we have for $x\in X$
\[
|T_{\alpha}(S_{\alpha}(x))|\leq |T_{\alpha}||S_{\alpha}(x)|\leq
|T_{\alpha}|(\varepsilon|x|)=|T_{\alpha}(\varepsilon|x|)|\leq
{\varepsilon}^2|x|,
\]
so that $|T_{\alpha}(S_{\alpha}(x))|\leq {\varepsilon}^2|x|$. Note
that by \cite{BK, CD}, for a central operator $T$ on a vector
lattice $X$, the modulus of $T$, $|T|$, exists and satisfies
$|T|(|x|)=|T(x)|$, for any $x\in X$. $\eproof$

\begin{proposition}\label{Z(X)2}
Suppose $(T_{\alpha})$ is a net of central operators on a
topological Riesz space $X$ converging (in the $\tau$-topology) to a
linear operator $T$. Then, $T$ is also central.
\end{proposition}
$\proof$ There is an $\alpha_0$ such that for each
$\alpha\geq\alpha_0$ and for each $x\in X$, we have
$|(T_{\alpha}-T)(x)|\leq |x|$. Fix an $\alpha\geq\alpha_0$. There
exists a positive real $\gamma_{\alpha}$ such that
$|T_{\alpha}(x)|\leq \gamma_{\alpha}|x|$. So, we have
\[
|T(x)|\leq |T_{\alpha}(x)|+|x|\leq
\gamma_{\alpha}|x|+|x|=(\gamma_{\alpha}+1)|x|.
\]
Therefore, $T$ is also a central operator. $\eproof$

\begin{proposition}\label{Z(X)3}
Let $X$ be a complete locally solid Riesz space. Then ${\sf Z}(X)$
is also complete with respect to the $\tau$-topology.
\end{proposition}
$\proof$ Let $(T_{\alpha})$ be a Cauchy net in ${\sf Z}(X)$ and $V$
be an arbitrary zero neighborhood in $X$. Fix  $x_0\in X$. Choose
$\delta>0$ such that $\delta|x_0|\in V$. There exists an $\alpha_0$
such that $|(T_{\alpha}-T_{\beta})(x)|\leq \delta |x|$ for each
$\alpha\ge\alpha_0$ and for each $\beta \ge \alpha_0$. So, we
conclude that $(T_{\alpha}-T_{\beta})(x_0)\in V$. This means that
$(T_{\alpha}(x_0))$ is a Cauchy net in $X$. Therefore, it is
convergent. Put $T(x)=\lim T_{\alpha}(x)$. Since this convergence
holds in ${\sf Z}(X)$, by Proposition \ref{Z(X)2}, $T$ is also
central. $\eproof$

Collecting the results of Propositions \ref{Z(X)1}, \ref{Z(X)2} and
\ref{Z(X)3}, we have the following

\begin{theorem}
Let $X$ be a complete locally solid Riesz space. Then, ${\sf Z}(X)$
is a complete unital topological algebra.
\end{theorem}

In addition, we have continuity of the lattice operations (defined
via formulas of \cite{BK, CD}) with respect to the assumed topology,
which is proved in the following theorem; in other words $({\sf
Z}(X), \tau)$ is a locally solid Riesz space.

\begin{theorem}
The lattice operations of ${\sf Z}(X)$ are continuous with respect
to the assumed topology.
\end{theorem}
$\proof$ By \cite{BK, CD}, the supremum and the infimum operations
in ${\sf Z}(X)$ are given by the formulas
\[
(T\vee S)(x)=T(x)\vee S(x) \mbox{ and } (T\wedge S)(x)= T(x)\wedge
S(x), \, \,  T,S\in {\sf Z}(X), x\in X^{+}.
\]
Let $(T_{\alpha})$ and $(S_{\alpha})$ be two nets of central
operators which are convergent to operators $T$ and $S$ in the
$\tau$-topology, respectively. Let $\varepsilon>0$ be arbitrary.
There are some $\alpha_0$ and $\alpha_1$ such that for each $x\in
X$, we have $|(T_{\alpha}-T)(x)|\leq \frac{\varepsilon}{2}|x|$ for
each $\alpha\geq\alpha_0$ and $|(S_{\alpha}-S)(x)|\leq
\frac{\varepsilon}{2}|x|$ for each $\alpha\geq\alpha_1$. Pick an
$\alpha_2$ with $\alpha_2\geq\alpha_0$ and $\alpha_2\geq\alpha_1$.
Then for each $\alpha\geq\alpha_2$, by using the Birkhoff's
inequality (for example, see \cite{AB}), for each $x\in X^{+} =
\{y\in X:y\ge 0\}$,
\[
|(T_{\alpha}\vee S_{\alpha})(x)-(T\vee S)(x)|=|(T_{\alpha}\vee
S_{\alpha})(x)-(T_{\alpha}\vee S)(x)+(T_{\alpha}\vee S)(x)-(T\vee
S)(x)|
\]
\[
\leq |T_{\alpha}(x)\vee S_{\alpha}(x)-T_{\alpha}(x)\vee
S(x)|+|T_{\alpha}(x)\vee S(x)-T(x)\vee S(x)|
\]
\[
\leq |(T_{\alpha}-T)(x)|+|(S_{\alpha}-S)(x)|\leq \varepsilon x.
\]
If $x$ is not positive, by using $x=x^{+}-x^{-}$, we get
\[
|(T_{\alpha}\vee S_{\alpha})(x)-(T\vee S)(x)|\leq 2 \varepsilon|x|
\]
for sufficiently large $\alpha$. Since, the lattice operations in a
Riesz space can be obtained via each other, we conclude that all of
them are continuous in ${\sf Z}(X)$ with respect to the
$\tau$-topology. $\eproof$

\medskip
\subsection{Order bounded below operators}

\medskip
In this subsection we present a new approach to order bounded below
operators on a topological Riesz space. The concept of a bounded
below operator on a Banach space has been studied widely (for
example, see Section 2.1 in \cite{AA}), and the norm of the operator
was used essentially in its definition. In a Riesz space we have the
concepts order and modulus which help us to have a different vision
to these operators. On the other hand, when we deal with a
topological Riesz space and an operator on it, the concept
"neighborhoods" enables us to consider a topological view for an
order bounded below operator. With the topology introduced in the
previous subsection we extend some results of \cite[Section 2.1]{AA}
to topological Riesz spaces.

A linear operator $T$ on a topological Riesz space $X$ is said to be
\emph{order bounded below} if $|T(x)|\ge \gamma |x|$ for some
positive real number $\gamma$ and each $x\in X$. The class of all
order bounded below operators on a topological Riesz space $X$ is
denoted by ${\sf O}(X)$. We consider the $\tau$-topology on it. It
is easy to see that every order bounded below operator is
one-to-one; also, every central order bounded below operator on a
locally solid Riesz space has a closed range.

In the following lemma we show that if a linear operator is
sufficiently close to an order bounded below operator with respect
to the $\tau$-topology, then it is also order bounded below. So, the
set of all central operators on a topological Riesz space $X$ which
are also order bounded below is an open subset of ${\sf Z}(X)$.

\begin{lemma}
If a linear operator $S$ on a topological Riesz space $X$ is
sufficiently close to an order bounded below operator $T$ with
respect to the $\tau$-topology, then $S$ is also order bounded
below.
\end{lemma}
$\proof$ Since $T$ is an order bounded below operator on $X$ there
is a positive real number $\gamma$ such that $|T(x)|\ge \gamma |x|$
for each $x\in X$. Choose a linear operator $S$ such that
$|(T-S)(x)|\leq \frac{\gamma}{2}|x|$. Then, for each $x\in X$, we
have
 \[
 |Sx|=|(S-T)(x)+T(x)|\ge \gamma |x|-\frac{\gamma}{2}|x|=\frac{\gamma}{2}|x|,
 \]
 so that $|S(x)|\ge \frac{\gamma}{2}|x|$, which completes the proof.
 $\eproof$

Now, we assert the main theorem of this subsection which is an
extension of Theorem 2.9 from \cite{AA} to a locally solid Riesz
space.

 \begin{theorem}
Let $X$ be a complete locally solid Riesz space and $T$ a continuous
order bounded below operator on $X$. If a net $(T_{\alpha})$ of
surjective operators converges in the $\tau$-topology to $T$, then
$T$ is also surjective.
\end{theorem}
$\proof$ Fix a constant $\gamma>0$ such that $|T(x)|\ge 2\gamma |x|$
for each $x\in X$. There is an $\alpha_0$ with
$|(T_{\alpha}-T)(x)|\leq \gamma|x|$ for each $x\in X$ and for each
$\alpha\ge\alpha_0$, so that $|T_{\alpha}(x)|\ge \gamma |x|$. Pick
$y\in X$. To show that $T$ is surjective, we have to show that there
exists $x\in X$ such that $T(x)=y$. There is a net $\{x_{\alpha}\}$
in $X$ such that $T_{\alpha}(x_{\alpha})=y$. From $\gamma
|x_{\alpha}|\leq |T_{\alpha}(x_{\alpha})|=|y|$, it follows that
$|x_{\alpha}|\leq \frac{1}{\gamma}|y|$.

We claim that $\{x_{\alpha}\}$ is a Cauchy net in $X$. Let $V$ be an
arbitrary zero neighborhood in $X$. Choose $\delta>0$ such that
$\delta \frac{1}{\gamma}|y|\in V$. Then for sufficiently large
$\alpha$ and $\beta$ we have
\[
2\gamma |x_{\alpha}-x_{\beta}|\leq |T(x_{\alpha})-T(x_{\beta})|\leq
|T(x_{\alpha})-y|+|y-T(x_{\beta})|
\]
\[
=|T(x_{\alpha})-T_{\alpha}(x_{\alpha})|+|T_{\beta}(x_{\beta})-T(x_{\beta})|
\leq \frac{\delta}{2}|x_{\alpha}|+\frac{\delta}{2}|x_{\beta}|\leq
\delta \frac{1}{\gamma}|y|.
\]
So, we conclude that $\{x_{\alpha}\}$ is really a Cauchy net in $X$.

Suppose that $(x_{\alpha})$ converges to $x$. Thus, it follows that
$(T(x_{\alpha}))$ converges to $T(x)$. On the other hand, from the
equality
\[
|T(x_{\alpha})-y|=|T(x_{\alpha})-T_{\alpha}(x_{\alpha})|\leq \delta
|x_{\alpha}|\leq \delta\frac{1}{\gamma}|y|
\]
it follows that $(T(x_{\alpha}))$ converges to $y \in X$, hence
$T(x)=y$. This complete the proof.
 $\eproof$

\footnotesize{
}

\end{document}